\documentclass[11pt,a4paper]{article}

\usepackage{bezier,amsfonts,amssymb,graphicx,amsthm,url}%5,mathtools}
\usepackage[english]{babel}
\def\lesta{ \hfill $\Box$ \bigskip}

\newcommand{\HH}{{\cal H}}

\newcommand{\uchi}{\overline{\chi}}

\newtheorem{theorem}{Theorem}[section]

\newtheorem{proposition}[theorem]{Proposition}

\newtheorem{lemma}[theorem]{Lemma}

\begin{document}

\title{Non-monochromatic non-rainbow colourings of $\sigma$-hypergraphs}
\author{Y. Caro \\ Deapartment of Mathematics\\ University of Haifa-Oranim \\ Israel \and J. Lauri \\ Department of Mathematics \\ University of Malta
\\ Malta
\\ \\ \emph{To the memory of Lucia Gionfriddo}}
\date{ }
\maketitle

\begin{abstract}
\begin{small}
One of the most interesting new developments in hypergraph colourings in these last few years has been Voloshin's notion of colourings of mixed hypergraphs. 
In this paper we shall study a specific instance of Voloshin's idea: a non-monochromatic non-rainbow (NMNR) colouring of a hypergraph is a colouring of its vertices such that every edge has at least two vertices coloured with different colours (non-monochromatic) and no edge has all of its vertices coloured with distinct colours (non-rainbow). Perhaps the most intriguing phenomenon of such colourings is that a hypergraph can have gaps in its NMNR chromatic spectrum, that is, for some $k_1 < k_2 < k_3$, the hypergraph is NMNR colourable with $k_1$ and with $k_3$ colours but not with $k_2$ colours.

Several beautiful examples have been constructed of NMNR colourings of hypergraphs exhibiting phenomena not seen in classical colourings. Many of these examples are either \emph{ad hoc} or else are based on designs. The latter are difficult to construct and they generally give uniform $r$-hypergraphs only for low values of $r$, generally $r=3$. In this paper we shall study the NMNR colourings of a type of $r$-uniform hypergraph which we call $\sigma$-hypergraphs. 
The attractive feature of these $\sigma$-hypergraphs is that they are easy to define, even for large $r$, and that, by suitable modifications of their parameters, they can give families of hypergraphs which are guaranteed to have NMNR spectra with gaps or NMNR spectra without gaps.  
These $\sigma$-hypergraphs also team up very well with the notion of colour-bounded hypergraphs recently introduced by Bujt{\'a}s and Tuza to give further control on the appearance of gaps and perhaps explain better the existence of gaps in the colouring of mixed hypergraphs.

\end{small}
\end{abstract}

\section{Introduction}

A \emph{hypergraph} $\HH$ is a finite set $V(\HH)$ of vertices together with a family $E(\HH)$ of subsets of $\HH$ called edges. When all the subsets are of the same size $r$ we say that $\HH$ is an \emph{$r$-uniform hypergraph} or simply a \emph{uniform hypergraph}. A \emph{$k$-non-monochromatic-non-rainbow} ($k$-NMNR) colouring of a uniform hypergraph $\HH$ is a colouring of the vertices of $\HH$ using exactly $k$ distinct colours such that no edge has all its vertices coloured with the same colour (a \emph{monochromatic edge})and no edge has all its vertices coloured with different colours (a \emph{rainbow edge}). 

These colourings are an important instance of what are sometimes called Voloshin colourings of mixed hypergraphs which is a generalisation of the classical colourings of hypergraphs \cite{voloshin02}. In Voloshin's notation, the edges of a \emph{mixed hypergraph} are defined to be of two types: the ${\cal D}$-edges and the ${\cal C}$-edges. In a colouring of a mixed hypergraph, the vertices of any $\cal D$-edge cannot all be given the same colour (that is, some two vertices must have $\cal D$ifferent colours), while the vertices in
a $\cal C$-edge cannot all be given different colours, (that is, some two vertices must have $\cal C$ommon colour). What we consider in this paper is the case of hypergraphs which, apart from being uniform, also have the property that every edge is both in class $\cal C$ and in class $\cal D$. Such hypergraphs are usually called bihypergraphs \cite{voloshin02,gionfriddo04}. 

The definition of mixed hypergraphs and much of the early theory itself of colourings of mixed hypergraphs has been introduced and developed by Vitaly Voloshin. The book \cite{voloshin02} gives an excellent picture of the situation up to 2002. The vitality of Voloshin's idea can be seen from the range of important papers published on the topic, such as the earlier ones, like
\cite{voloshin95,voloshinurl,colburn&99,kundgen&00,dvorak&kra01,gionfriddo&02,gionfriddo&02}, or more recent examples
\cite{kral&03,kral04,gionfriddo04,jiang&05,gyori&08} and now even applications  
\cite{tuza13,zhang&13,jaffe&12}. The regularly updated web site \url{http://http://spectrum.troy.edu/voloshin/publishe.html} is an ideal source for up-to-date-material on colourings of mixed hypergraphs.

Unlike classical colourings of hypergraphs, it is possible for a hypergraph not to have a $k$-NMNR colouring for any $k$. Two examples would be hypergraphs all of whose edges have size 2 or complete $r$-uniform hypergraphs on $n>(r-1)^2$ vertices. We shall take these to be trivially non-NMNR-colourable hypergraphs---some non-trivial cases will be considered below. 

But, if $\HH$ does have a $k$-NMNR colouring for some $k$ then, as in classical colourings of hypergraphs, there is a number $\chi=\chi(\HH)$ such that $\HH$ has a $\chi$-NMNR colouring but it does not have a $k$-NMNR colouring for any $k<\chi$. However, unlike classical colourings, there could be a number $k <|V(\HH)|$ such that $\HH$ does not have a 
$k'$-colouring for any $k'>k$. The largest such number $k$ for which $\HH$ does have a $k$-NMNR colouring is denoted by $\uchi=\uchi(\HH)$. We call $\chi$ and $\uchi$ the \emph{lower} and the \emph{upper} NMNR chromatic numbers of $\HH, respectively$. 

A property of NMNR colourings which we find very intriguing is that there could be a number $k$ with $\chi<k<\uchi$ such that $\HH$ is not $k$-NMNR colourable. 
More specifically, 
we say (differring slightly from the notation in \cite{voloshin02}) that  the \emph{NMNR spectrum} of a hypergraph $\HH$ is the sequence, in increasing order, of all $k$ such that $\HH$ has a $k$-NMNR colouring. The first and last terms of this sequence are clearly $\chi=\chi(\HH)$ and $\uchi=\uchi(\HH)$, respectively. If the spectrum contains all integers between $\chi$ and $\uchi$ we then say that it is \emph{continuous}. Otherwise we say that $\HH$ has a \emph{gap} in its NMNR spectrum or has a \emph{broken} NMNR spectrum. 
The theme of gaps in the NMNR spectrum will be the main topic of study in this paper.

Unless otherwise stated, all hypergraphs $\HH$ considered will be uniform. Any graph theory terms not defined here can be found in \cite{bondy&mur08} and any undefined terms on hypergraphs can be found in \cite{voloshin02}.        

\section{$\sigma$-hypergraphs}

We now come to the definition of the main protagonists of this paper. There have been various instances in the literature where researchers have tried to generalise the idea of mixed colourings, and we shall be meeting an important one (colour-bounded hypergraphs) at the end of this paper. Perhaps the most ambitious attempt at defining hypergraphs on which a full generalisation of mixed colourings can be investigated is the idea of ``pattern hypergraphs'' as defined in \cite{dvorak&10}. What we shall instead do here is present a simple model of a hypergraph which is very easy to define and which can give us different NMNR colouring  properties by a simple adjustment of its parameters. These hypergraphs also have the nice property of integrating in a very interesting way with colour-bounded hypergraphs as we shall eventually see.  

So, let $n,r,q$ be positive integers, and let $\sigma$ be a partition of the integer $r$. The $\sigma$-hypergraph $\HH=\HH(n,r,q|\sigma)$ is defined as follows. It has $nq$ vertices and these are partitioned into $n$ subsets called \emph{classes} each containing $q$ vertices. If the classes are denoted by $V_1,V_2,\ldots,V_n$ then a subset $R$ of $V(\HH)$ of size $r$ is an edge if the partition of $r$ formed by the non-zero cardinalities $|R\cap V_i|, 1\leq i\leq n$ is $\sigma$. We shall assume, throughout, that $r\geq3$.

We shall denote by $\Delta(\sigma), \delta(\sigma)$ and $s(\sigma)$ the size of the largest part of the partition $\sigma$, the size of the smallest part of $\sigma$, and the number of parts (summands) of $\sigma$, respectively. For example, if $r=15$ and $\sigma=(1,1,2,3,3,5)$, 
every edge $K$ will be a 15-subset of the vertices of the hypergraph which intersects exactly six of the classes giving intersections of sizes 1, 1, 2, 3, 3 and 5, respectively. In this case, $\Delta(\sigma)=5$, $\delta(\sigma)=1$ and $s(\sigma)=6$.

We shall always assume that $s(\sigma)\geq2$ so that we avoid the possibility of an edge lying completely within one class. 
To simplify notation, when $\sigma$ is the partition $(1^{r-p},p)$ we denote $\HH(n,r,q|\sigma)$ by $\HH(n,r,q,p)$. Therefore the hypergraph above,
$\HH(2n,r)$, can also be denoted as the $\sigma$-hypergraph $\HH(n,r,2,2)$

\medskip\noindent
We have found that $\sigma$-hypergraphs are very versatile indeed. They can give us families of hypergraphs which exhibit a variety of interesting colouring properties of mixed hypergraphs. For example, it is easy to construct families of hypergraphs which do not have any NMNR colouring. Recall that any edge of $\HH(n,r,q,1)$ is any $r$-subset of vertices taken from $r$ different classes.

\begin{proposition}
Consider the $\sigma$-hypergraph $\HH=\HH(n,r,q,1)$. 
If $n\geq (r-1)^2 +1$ then there is no $k$-NMNR colouring for any value of $k$.
\end{proposition}

\noindent
{\bf Proof.} This is a simple application of the pigeonhole principle. Suppose that the vertices of $\HH$ have been coloured using the $k$ colours $\{1,2,\ldots,k\}$, for any $k\geq1$. Take a subset of vertices $\{v_1,\ldots,v_n\}$ each one from a different class of $\HH$. Place them into subsets $S_1,...,S_k$ where each $S_i$ conatins all those vertices from $\{v_1,\ldots,v_n\}$ coloured with colour $i$. Suppose $\HH$ contains no monochromatic edge. Therefore each subset $S_i$ has size at most $r-1$. But $n\geq (r-1)^2+1$, therefore there are at least $r$ non-empty subsets $S_i$. But then, $r$ vertices taken from different non-empty subsets $S_i$ give a rainbow edge. Therefore $\HH$ is not $k$-NMNR colourable.  
\lesta

The NMNR colouring properties of $\HH(n,r,q,1)$ for $n\leq (r-1)^2$ are not difficult to analyse. Since our aim in this paper is to use $\sigma$-hypergraphs in order to study gaps in the NMNR spectrum, we merely invite the interested reader to look at this situation.

\section{Gaps in the NMNR spectrum of $\sigma$-hypergraphs}

In this section, using $\sigma$-hypergraphs we shall be constructing families of hypergraphs with various properties related to gaps in the NMNR spectrum. But first we shall give a simple example of this counter-intuitive phenomenon of gaps. This example will motivate our general construction. Following this, we shall show that, in general, $\sigma$-hypergraphs with $\Delta(\sigma)>\delta(\sigma)=1$ can have broken NMNR spectrum. We shall then give a detailed analysis of such a family of hypergraphs. Finally, we shall show that a $\sigma$-hypergraph with $\delta(\sigma)\geq 2$ can never have a broken spectrum.

\subsection{A simple example} \label{sec:simpleExample}

Let us consider the following family of hypergraphs which we denote by $\HH(2n,3)$. It has $2n$ vertices, $n\geq4$, which are partitioned into $n$ pairs. A edge is defined to be a set of three vertices such that exacly two are in the same pair. Clearly, $\HH(2n,3)$ is 2-NMNR colourable: just colour the two vertices in each pair with the colours 1 and 2. Also, $\HH(2n,3)$ is clearly $n$-NMNR colourable: just give the two vertices in each pair the same colour, using a different colour for each pair. However, $\HH(2n,3)$ is not $k$-NMNR colourable for any $3\leq k \leq n-1$. Because suppose it is. If two vertices in the same pair are given two different colours, say 1 and 2, then find a vertex with colour 3. These three vertices give a rainbow edge. On the other hand, if every two vertices in a pair are given the same colour then, by the pigeon-hole prionciple, there must be two pairs such that all their vertices are given the same colour. Taking two vertices from one of these pairs and the third from the other gives a monochromatic edge. 

In fact, the full spectrum of $\HH(2n,3)$ is $\{2,n\}$, because if we use more than $n$ colours, then certainly one of the $n$ pairs contains distinct colours, say 1 and 2. Since we are using at least $n+1\geq3$ colours, some other pair contains a vertex coloured 3. But then these three vertices form a rainbow edge.

Our $\sigma$-hypergraphs are a generalisation of this simple example.
The hypergraph $\HH(2n,3)$, which we have just defined, is the $\sigma$-hypergraph $\HH(n,r,q|\sigma)$ with $r=3$, $q=2$ and $\sigma=(1,2)$.  This particular class of $\sigma$-hypergraphs can be extended as follows: let there still be $n$ classes of two vertices each, but let the edges be all $r$-subsets which contain exactly two vertices coming from the same class. This therefore becomes the $\sigma$-hypergraph $\HH(n,r,2|\sigma=(1^{r-2},2))$, which we denote, for short, by $\HH(2n,r)$. We shall see below that this is a very interesting family of hypergraphs when we shall be discussing in detail its NMNR spectrum. 

\subsection{The monochromatic zone}

We shall first note that, for any $\sigma$-hypergraph with $\Delta(\sigma)>1$ there is a range of consecutive integers $k$ for which there is always a $k$-NMNR colouring. We shall refer to this range of integers as the monochromatic zone.

So, consider a $\sigma$-hypergraph $\HH(n,r,q|\sigma)$, always with $s=s(\sigma)\geq2$. We also assume here that $\Delta(\sigma)>1$. Colour all the vertices in each of the $n$ classes of vertices with the same colour. How many colours can be used this way to give us an NMNR colouring? Suppose first that all classes are given different colours. Then, since no edge is contained in a single class, there can be no monochromatic edge. Also, since not all parts of the partition are equal to 1, because $\Delta(\sigma)>1$, there can be no rainbow edge. Therefore we have an $n$-NMNR colouring. At the other end, suppose, in this scheme, that no colour repeats in more than $s-1$ classes. Then again, we have no monochromatic or rainbow edges, therefore giving us a $\lceil \frac{n}{s-1} \rceil$-NMNR colouring. Clearly, this scheme gives us a $k$-NMNR colouring for all $\lceil \frac{n}{s-1} \rceil \leq k \leq n$ and we call this interval the \emph{monochromatic zone}, which is therefore a gap-free region within the NMNR spectrum. Since we shall be using this result often, we record it here as a proposition.

\begin{proposition} \label{prop:monochromaticZone}
Let $\HH$ be a $\sigma$-hypergraph with $\Delta(\sigma)>1$ and $s(\sigma)\geq2$. Then $\HH$ is $k$-NMNR colourable for all $k$ in the interval 
\[\left\lceil \frac{n}{s-1} \right\rceil \leq k \leq n.\]
\end{proposition}

\subsection{Gaps in the case $\Delta(\sigma)>\delta(\sigma)=1$}

We shall now show, by means of an example, that when $\Delta(\sigma)>\delta(\sigma)=1$
a certain choice of the parameter $q$ gives $\sigma$-hypergraphs with a gap in the spectrum.
We recall that we write $s=s(\sigma)$ for the number of parts of $\sigma$ and that we always assume that $s\geq2$. We shall also write $\Delta$ for $\Delta(\sigma)$.

\begin{theorem}
Let $n$ be greater than $r(s-1)+s$, and let $q=(r-1)(\Delta-1)$. Let $\sigma$ be a partition of $r$ with $\delta(\sigma)=1$, and let $\HH=\HH(n,r,q|\sigma)$. Then $\HH$ has a gap in its spectrum between $r-1$ and the monochromatic zone. More specifically, $\HH$ is 
\begin{enumerate}
\item not $k$-NMNR colourable for $k\leq r-2$;
\item $(r-1)$-NMNR colourable;
\item not $r$-NMNR colourable;
\item $k$-NMNR colourable in the monochromatic zone $\lceil n/(s-1)\rceil \leq k \leq n$.
\end{enumerate}
\end{theorem}

\noindent
{\bf Proof.} {\it Case 1: $k\leq r-2$.}\newline
So, suppose first that $\HH$ has been coloured with $k\leq r-2$ colours. Then, in every class, there is a colour which appears at least $\Delta$ times. This is because, if not, then the number of vertices in some class would be at most $(\Delta-1)k \leq (\Delta-1)(r-2) < (\Delta-1)(r-1) = q$, where $q$ is the number of vertices in a class.
 
Hence, since $n\geq r(s-1) + s > k(s-1)+1$, for $k\leq r-2$, it follows by the pigeon-hole principle, that one colour appears at least $\Delta$ times in at least $s=s(\sigma)$ classes. But this means that there is a monochromatic edge. Therefore $\HH$ is not $k$-NMNR colourable for $n\geq k(s-1) + 1$ and $k\leq r-2$, in particular for $n \geq (r-2)(s-1)+1$.

\medskip\noindent
{\it Case 2: $k=r-1$.}\newline
Next we show that $\HH$ is $(r-1)$-colourable. Colour each class as 
\[(1^{\Delta-1}, 2^{\Delta-1}, \ldots, (r-1)^{\Delta-1}),\] 
that is, each colour $1,2,\ldots,r-1$ is repeated exactly $\Delta-1$ times in each class. Then there can be no monochromatic edge since no colour appears more than $r-1$ times in any class. Also, since only $r-1$ colours are used, there cannot be any rainbow edge.

\medskip\noindent
{\it Case 3: $k=r$.}\newline We now have to prove now is that $\HH$ is not $r$-NMNR colourable. For the rest of the proof, let the parts of the partition $\sigma$ be denoted by $\Delta=a_1, a_2,\ldots, a_s=1$ in non-increasing order. Observe that $a_1+\ldots+a_{s-1} = r-1$, a fact that will be needed later.

So, suppose we have the $r$ colours $1,2\ldots,r$. Observe first that the same colour can appear at least $\Delta$ times in at most $s-1$ classes, otherwise we obtain a monochromatic edge. Also, since $r$ colours are used, there are at most $r(s-1)$ classes which contain some colour repeated at least $\Delta$ times.

But the number of classes $n$ is at least $r(s-1) +s$. Therefore there are at least $s$ classes in which no colour repeats more than $\Delta-1$ times. Let us call such classes ``colourful''. Since $q$, the size of each class, equals $(r-1)(\Delta-1)$, it follows that each colourful class contains either $r-1$ or $r$  distinct colours, which is why we call them colourful. We shall assume that these $s$ colourful classes are listed as $V_1, V_2, \ldots, V_s$ in some arbitrary order. 

We now consider two cases. Firstly, suppose some colour $x$ is missing from all the colourful classes. Therefore $x$ must appear at least once in some other class. We may assume, by re-numbering, that $x=r$, and that therefore each colourful class contains the colours $1,2,\ldots,r-1$. In that case we can form a rainbow edge $K$ as follows: $K$ contains $\Delta=a_1$ distinct colours from $V_1$, $a_2$ new distinct colours from $V_2$, all different from the ones previously chosen, and so on up to $a_{s-1}$ new distinct colours from the $V_{s-1}$.
Note that this process of choosing greedily the distinct colours class after class is possible due to the above observation that $a_1+\ldots+a_{s-1}=r-1$. 

Finally, we take the colour $r$ from  whichever class it appears in (which is not colourful). And this gives us the rainbow edge $K$, which makes this case impossible.

We now consider the second case, that is, when each of the colours $1,2,\ldots,r$ appears at least once in some colourful class. Again, we shall construct a rainbow edge $K$. The first $r-1$ distinct colours in $K$ are chosen by the previous greedy fashion: choose $a_j$ new distinct colours from the $V_j$, for $j=1,2,\ldots,s-1$. We just need to assign the last colour to $K$.

Note that $K$ already contains all colours except one, call it $x$, say. If $x$ appears in $V_s$ then we assign it as the last colour to $K$. Otherwise, $x$ must be in some $V_j$ from which we have chosen $a_j$ distinct colours. But then these $a_j$ colours must all appear in $V_s$ since $V_s$ misses colour $x$ and it can miss at most one colour. Therefore we re-assign to $K$ the colour $x$ from $V_j$ instead of some colour $y$ from the $a_j$ colours previously assigned, and then we assign to $K$ the colour $y$ from $V_s$. This, again, gives us the rainbow edge $K$, which is the final contradiction we required.

\medskip\noindent
{\it Case 4: $\lceil n/(s-1) \leq k \leq n$.}\newline
We already know that $\HH$ is $k$-NMNR colourable for $k$ in the monochromatic zone, so we have nothing to probe here.
\lesta

\subsection{The hypergraph $\HH(2n,r)$: full calculation of a spectrum with a gap}

In the previous subsection we have simply shown that a particular family of $\sigma$-hypergraphs does have a gap but we did not try to determine the exact extent of the gap. In this section we shall study the interesting and relatively simple $\sigma$-hypergraph $\HH(2n,r)$, which is an extension of the elementray example we started with in order to introduce the notion of a broken NMNR spectrum, and we shall show that we can determine exactly its NMNR spectrum.

Recall that $\HH(2n,r)=\HH(n,r,2|(1^{r-2},2))$, that is, it is the hypergraph on $2n$ vertices which are partitioned into $n$ classes of size 2, and its edges are all $r$-subsets which contain two vertices from the same class and all the other $r-2$ vertices from distinct classes.

We shall first deal with the case when $r\geq6$, following which we shall briefly treat the details which arise for smaller values of $r$.

\begin{theorem} \label{thm:theMonster}
Let $\HH=\HH(2n,r)$ with $r\geq6$. Suppose also that $n\geq 2(r-2)(r-1)+1$.
Then,
\begin{enumerate}
\item $\HH$ is not $k$-NMNR colourable for $k>n$;
\item $\HH$ is $k$-NMNR colourable for $\lceil n/(r-2) \rceil \leq k \leq n$;
\item $\HH$ is not $k$-NMNR colourable for $2r-5 \leq k < \lceil n/(r-2) \rceil$;
\item $\HH$ is $k$-NMNR colourable for $2\leq k \leq 2r-6$; 
\end{enumerate}
\end{theorem}

\noindent
{\bf Proof.} We let the classes containing the vertices of $\HH$ be $V_1,V_2,\ldots,V_n$ in some arbitrary order. Recall that each $V_i$ has size equal to 2.

\medskip\noindent
{\it Case 1: $k>n$} \newline
We shall first show that, if the number of colours $k$ used is more than $n$, then there will always be a rainbow edge. So, suppose $k>n$ and $\HH$ is coloured using $k$ colours. Therefore there is at least one class which contains two different colours. We may assume, without loss of generality, that $V_1$ contains the colours $1,2$. We construct a rainbow edge $K$ as follows. First let $K$ contain the two vertices from $V_1$ and another vertex coloured 3. We may assume, by renumbering, that this vertex is from $V_2$. We keep choosing a vertex with a new colour from some new class until we obtain a rainbow $r$-edge. How do we know that this is certainly possible? Suppose not. Therefore suppose that we have reached a stage where $K$ conatains $j$ differently coloured vertices from the sets $V_1,V_2,\ldots,V_{j-1}$ with $j<r$, and that we cannot find amongst the classes $V_j,\ldots,V_n$ any vertex with a different colour from those already in $K$. But then, the number of colours used in the colouring is at most $2(j-1)$. However,
$2(j-1) < 2(r-1) < (2r-4)(r-1) \leq n < k$, 
which is a contradiction. Therefore we can always find a rainbow edge, which shows that $\HH$ is not $k$-NMNR colourable for $k>n$.

\medskip\noindent 
{\it Case 2: $\lceil n/(r-2) \rceil \leq k \leq n$} \newline
We already know that $\HH$ is $k$-NMNR colourable for $\lceil n/(r-2) \rceil \leq k \leq n$ since this is the monochromatic zone, so we have nothing to prove here. 

\medskip\noindent
{\it Case 3: $2r-5 \leq k < \lceil n/(r-2) \rceil$} \newline
So we now have to show that $\HH$ is not $k$-NMNR colourable for $2r-5 \leq k < \lceil n/(r-2) \rceil$. We shall first tackle the interval $2r-3 \leq k < \lceil n/(r-2) \rceil$. The cases $k=2r-4$ and $k=2r-5$ will then be treated separately below.

Firstly, if all the classes $V_i$ were monochromatic, then, since no more than $r-2$ classes can be monochromatic of the same colour (otherwise we get a monochromatic edge), the number of colours used would be at least $\lceil n/(r-2) \rceil$. But this is impossible since $k < \lceil n/(r-2) \rceil$. So there is some class, say $V_1$, which is not monochromatic.

But we have $k\geq 2r-3$ colours. So it is  not possible that the remaining $2r-5$ colours (the ones not in $V_1$) are contained in only $r-3$ classes, which between them would contain at most $2r-6$ vertices. Therefore there are $r-2$ classes which, between them, certainly contain $r-2$ distinct colours, one colour in each class. But then these classes, together with $V_1$, give a rainbow edge. So $\HH$ is not $k$-NMNR colourable for $2r-3 \leq k < \lceil n/(r-2)\rceil$.

\smallskip
We now proceed to treat the remaining cases, $k=2r-5$ and $k=2r-4$. First we show that $\HH$ is not $(2r-4)$-NMNR colourable. Suppose the vertices of $\HH$ have been coloured with $2r-4$ colours. First we shall show that there must be a non-monochromatic class. For, suppose all classes are monochromatic. Since there are $n$ classes and $n$ is at least $(2r-4)(r-1)+1$, some colour must appear in at least $r-1$ monochromatic classes. But then we have a monochromatic edge.

So, suppose that $V_1$ is a non-monochromatic class of $\HH$, and suppose, without loss of generality, that the colours appearing in $V_1$ are 1 and 2.  We shall construct a rainbow edge in the now familiar manner.  We start with the two vertices from $V_1$. Then we look for another colour from a different class; we may assume, by re-numbering colours and classes, that we find a vertex with colour 3 from $V_2$. In the next step we find, say, a vertex coloured 4 from $V_3$, and so on.  If we can do this up to colour $r$ then we have a rainbow edge.  So suppose we only manage to arrive up to some colour $j\leq r-1$ from  $V_{j-1}$.  Since we cannot proceed any further, no other new colour is available because it is found only in the classes $V_1$ to $V_{j-1}$. Therefore these classes between them contain all $2r-4$ colours, and therefore $j=r-1$. We therefore have that each class $V_i$, 
for $i\in \{1,\ldots,r-2\}$, has two vertices coloured  differently and all colours found in these classes are diffferent. 

But then we can again get a rainbow edge as follows.  Start with any vertex coloured $x$, say, from some class other than $V_1,\ldots,V_{r-2}$. Then, from each of the classes $V_1$ to $V_{r-2}$ pick one vertex, being careful not to pick the only vertex in these classes coloured $x$;  and then choose one other vertex from these same first $r-2$ classes (thereby repeating one of them once), again not choosing the one vertex coloured $x$.  This gives a rainbow edge. Therefore $\HH$ is not $(2r-4)$-colourable.

\smallskip
Similarly, we now show that $\HH$ is not $(2r-5)$-NMNR colourable. Again, if all classes were monochromatic, then, since there are $n$ classes and $n$ is at least $(2r-5)(r-1)+1$, some colour must appear in at least $r-1$ monochromatic classes giving a monochromatic edge. Therefore we may assume that $\HH$ has a non-monochromatic class $V_1$ with vertices coloured 1 and 2, respectively. As above, we start with the two vertices in $V_1$ and then construct a rainbow edge by choosing a sequence of vertices, each with a new colour and from a new class and we assume that we stop short of a rainbow edge when we have chosen a colour $j$ from $V_{j-1}$ for some $3\leq j <r$. Again, the classes $V_1$ to $V_{j-1}$ must therefore contain all the $2r-5$ colours used in the colouring, therefore $j=r-1$ and all these classes contain two differently coloured vertices,  except possibly for one of them which might be monochromatic. If all these classes are non-monochromatic, we proceed as before. So suppose there is a monochromatic class from amongst $V_1,\ldots,V_{j-1}$ which contains two vertices coloured $y$. If all the vertices in the classes $V_{r-1}$ to $V_n$ are coloured $y$ then $\HH$ contains a monochromatic edge. So, there is some vertex in these classes coloured $x\not=y$. Now, as before, obtain a rainbow edge starting with this vertex coloured $x$ and then using only vertices from the classes $V_1$ to $V_{j-1}$. Therefore $\HH$ does not have a $(2r-5)$-NMNR colouring.

\medskip\noindent
{\it Case 4: $2\leq k \leq 2r-6$} \newline
Finally, we have to show that there is always a $k$-NMNR colouring for $2 \leq k \leq 2r-6$. Let $t$ be a number such that $0 \leq t \leq r-4$. Divide the classes in $\HH$ into two categories: $V_1,\ldots, V_t$ and $V_{t+1},\ldots,V_n$ (if $t=0$ we assume that all the classes $V_1$ to $V_n$ are considered to be in the second category). We now colour the vertices of $\HH$ as follows. Each class from $V_{t+1},\ldots,V_n$ is given the two colours $1,2$. The vertices in the other classes are given distinct colours, all diferent from $1$ or $2$. Therefore this is a colouring which uses between 2 and $2 + 2(r-4) = 2r-6$ colours, depending on the value of $t$. With this colouring, an edge $K$ cannot be monochromatic. Also, $K$ can have at most $r-3$ vertices from the classes $V_1,\ldots,V_t$ (the maximum is attained if $K$ includes a whole class from this list). Therefore at least three other vertices must share the colours 1 and 2 from the classes $V_{t+1}$ to $V_n$, hence $K$ is not rainbow. We have therefore shown that $\HH$ is $k$-NMNR colourable for $2 \leq k \leq 2r-6$. 
\lesta

We are therefore left with the cases $r=3,4,5$ which need some special treatment not covered by Theorem \ref{thm:theMonster}. We shall consider these cases now.

\begin{proposition}
Consider the hypergraph $\HH(2n,r)$ for $r=3,4,5$ and $n\geq 2(r-2)(r-1)+1$. Then,
\begin{enumerate}
\item $\HH(2n,3)$ is $k$-NMNR colourable only for $k=2$ and $k=n$.
\item $\HH(2n,4)$ is $k$-NMNR colourable only for $k=2$, $k=3$ and for $k$ in the monochromatic zone, that is, $\lceil n/2 \rceil \leq k \leq n$.
\item $\HH(2n,5)$ is $k$-NMNR colourable only for $k=2,3,4$ and for $k$ in the monochromatic zone, that is, $\lceil n/3 \rceil \leq k \leq n$.
\end{enumerate}
\end{proposition}

\noindent
{\bf Proof.}
We have already shown in Section \ref{sec:simpleExample} that $\HH(2n,3)$ is $k$-NMNR colourable only for $k=2$ and $k=n$.

\medskip\noindent
Unlike the general case of Theorem \ref{thm:theMonster}, $\HH(2n,4)$ is $k$-NMNR colourable for $r=2r-5=3$. The colouring with three colours for $\HH(2n,4)$ is obtained by colouring all classes using the colours 1,2 except one class whose two vertices are both coloured 3. 

The proof that $\HH(2n,4)$ is not $k$-NMNR colourable for $4 \leq k < \lceil n/2 \rceil$ is as  that in Theorem \ref{thm:theMonster} but it requires a special treatment for $k=4$ (which corresponds to the case $2r-4$ in Theorem \ref{thm:theMonster}). 

Thus, suppose $\HH(2n,4)$ is coloured using four colours.  No colour can appear in three monochromatic classes, otherwise we would have a monochromatic edge. Therefore we can have at most eight monochromatic classes. But $n\geq 2(r-2)(r-1)+1>8$. Therefore there must be a non-monochromatic class, and we may take it to be $V_1$ whose vertices are coloured 1 and 2, respectively. Now, there is some vertex coloured differently in some other class. By renumbering colours and classes, we may assume that a vertex coloured 3 is in class $V_2$. If the colour 4 is found in some third class, then we have a rainbow edge. So we may suppose that colour 4 appears together with colour 3 in class $V_2$ and nowhere else. But now, take a vertex in $V_3$ coloured $x\in\{1,2,3,4\}$. We may assume, without loss of generality, that $x=1$. Then this vertex, together with the vertex coloured 2 in $V_1$ and the two vertices in $V_2$ give a rainbow edge. Therefore $\HH(2n,4)$ is not 4-NMNR colourable.

\medskip\noindent
Finally, the result for the hypergraph $\HH(2n,5)$ tells us that Theorem \ref{thm:theMonster} holds here. But again, the colourings for small values of $k$ need separate explaining. The $3$-NMNR coloring is once more obtained by colouring all classes with the colours 1,2 except for one monochromatic class whose vertices are given the colour 3. The $4$-NMNR colouring is obtained by colouring all classes except two with the colours 1,2. The remaining two classes are monochromatic: the vertices of one are both coloured 3 while the vertices of the other are both coloured 4. 

The proof that $\HH(2n,5)$ is not $k$-NMNR colourable for $5\leq k < \lceil n/3 \rceil$ is again similar to Theorem \ref{thm:theMonster} but it requires special treatment for the cases $k=5,6$ (corresponding to the cases $2r-5$ and $2r-4$, respectively, in Theorem \ref{thm:theMonster}).

So, suppose first that $\HH(2n,5)$ is coloured using six colours. No colour can appear in four monochromatic classes, otherwise we would have a monochromatic edge. Therefore we can have at most fifteen monochromatic classes. But $n\geq 2(r-2)(r-1)+1>15$. Therefore there must be a non-monochromatic class, and we may take it to be $V_1$ whose vertices are coloured 1 and 2, respectively. We may assume, by renumbering, that there is a vertex coloured 3 in $V_2$. We proceed this way, picking a vertex coloured with a new colour coming from a new class. If we can do this up to colour 5, then we have a rainbow edge. So the procedure must stop with colour $j<5$ in class $V_{j-1}$ (as usual, by renumbering of colours and classes). Since we cannot proceed further it follows that all colours greater than $j$ are to be found only in the classes $V_1$ to $V_{j-1}$. But six colours must be used, therefore, for the classes $V_1$ to $V_{j-1}$ to accommodate six colours,  $j$ must be equal to 4. We may assume, again by renumbering, that the vertices of $V_1, V_2, V_3$ are coloured, respectively, $(1,2)$, $(3,4)$ and $(5,6)$. But then, take any vertex $v$ from $V_4$ coloured $x$, say. We may assume, without loss of generality, that $x=1$. Therefore, taking the vertex $v$, the vertex coloured 2 from $V_1$, the two vertices from $V_2$ and any vertex from $V_3$ give a rainbow 5-edge. So $\HH(2n,5)$ is not 6-NMNR colourable.

Finally, we show that $\HH(2n,5)$ is not 5-NMNR colourable. The proof is similar but at one point requires some more care. Suppose $\HH(2n,5)$ has been coloured using five colours. As usual, we start by noting that there must be a non-monochromatic class which we take to be $V_1$ with two vertices coloured 1 and 2, respectively. As above, we begin constucting a rainbow edge by picking a sequence of vertices each with a new colour and coming from a new class. If we can do this up to colour 5 we get a rainbow edge. So the procedure must stop with colour $j<5$ in class $V_{j-1}$. Since we cannot proceed further it follows that all colours greater than $j$ are to be found only in the classes $V_1$ to $V_{j-1}$. But five colours must be used, therefore, for the classes $V_1$ to $V_{j-1}$ to accommodate five colours,  $j$ must be equal to 4. We now have to consider two cases. One case occurs when one of the classes $V_1, V_2, V_3$ is monochromtic, and the other when none of them is monochromatic. 

Let us start with the first case. We may assume, without loss of generality, that the vertices of $V_1, V_2, V_3$ are coloured, respectively, $(1,2)$, $(3,5)$ and $(4,4)$. It is clear that not all of the vertices in the classes $V_4$ to $V_n$ can be coloured 4, otherwise we would have a monochromatic edge. So suppose that some vertex $v$ in these classes is coloured $x\not=4$. We may assume, again without loos of generality, that $x=1$. But then, choosing the vertex $v$, the vertex coloured 2 in $V_1$, the two vertices in $V_2$, and one vertex from $V_3$, gives a rainbow 5-edge.

If, on the other hand, none of $V_1, V_2, V_3$ is monochromatic, we may assume that their vertices are coloured, respectively, $(1,2)$, $(3,4)$ and $(4,5)$. But then, let $v$ be a vertex coloured $x$ from $V_4$. It is easily checked that, for any $x$, one can find four appropriate vertices from $V_1, V_2, V_3$ which, together with $v$, form a rainbow 5-edge.

Therefore in any case, $\HH(2n,5)$ is not 5-NMNR colourable.
\lesta

\subsection{The case $\delta(\sigma)\geq2$: no gaps!}

We shall now prove the general result that, when $\delta(\sigma)\geq2$, a $\sigma$-hypergraph cannot have gaps in its NMNR spectrum.
We shall start with two technical lemmas which show how, under certain conditions, we can obtain a new NMNR colouring by re-colouring a given NMNR colouring. 

\begin{lemma} \label{lem:firstSwitchOfColours}
Let $\HH=\HH(n,r,q|\sigma)$ with $\delta(\sigma)\geq2$ and suppose we are given an NMNR colouring of $\HH$. Suppose $V$ is a class of $\HH$ and that all colours of $V$ are changed into a new colour $z$. Then the new colouring of $\HH$ is also an NMNR colouring. 
\end{lemma}

\noindent
{\bf Proof.} Any edge of $\HH$ which does not intersect $V$ remains a NMNR edge in the new colouring. Let $K$ be a edge which intersects $V$. Then, since $s(\sigma)\geq2$, $K$ intersects some other class $V'$. Therefore $K$ is non-monochromatic since it contains the colour $z$ which is not found in $V'$. Also, since $\delta(\sigma)\geq2$, $K$ contains the colour $z$ at least twice therefore it is non-rainbow. \lesta

\begin{lemma} \label{lem:secondSwitchOfColours}
Let $\HH=\HH(n,r,q|\sigma)$ with $\delta(\sigma)\geq2$ and suppose we are given an NMNR colouring of $\HH$. Suppose $V$ is a class of $\HH$ and that $x,y$ are two different colours in $V$ which do not appear in any other class of $\HH$. Suppose that all occurrences of the colours $x$ and $y$ in $V$ are changed into a new colour $z$. Then the new colouring of $\HH$ is also an NMNR colouring. 
\end{lemma}

\noindent
{\bf Proof.} As usual, if $K$ is a edge which does not interesect $V$, then it remains unchanged and therefore a NMNR edge in the new colouring. Therefore suppose that $K$ interesects $V$. If it did not contain any of the colours $x,y$ in the original colouring, it would again remain unchanged and therefore a NMNR edge in the new colouring. Therefore suppose first that $K$ contains at least two vertices coloured $x$ or $y$ in the original colouring. Therefore $K$ contains, in the new colouring, at least two vertices which are coloured $z$; therefore $K$ is non-rainbow. Also, since $z$ does not appear in any other class and $s(\sigma)\geq2$, $K$ contains some other colour from another class which is not $z$. Therefore $K$ is non-monochromatic.

Now suppose that $K$ contains only one vertex which is coloured $z$ in the new colouring. Then, as before, $K$ is non-monochromatic, since $s(\sigma)\geq2$ and $z$ is a new colour appearing only in $V$. Also, suppose, without loss of generality, that the vertex coloured $z$ in $K$ was coloured $x$ in the original colouring---therefore only one vertex was coloured $x$ in $K$. Recall that $K$ was non-rainbow, therefore it had at least two vertices coloured some colour $c$, with $c$ not equal to $x$ since $x$ does not appear in any class other than $V$. But this means that in the new colouring, these two vertices are still coloured $c$, therefore $K$ is still non-rainbow.   
\lesta

\begin{proposition} \label{th:noUpperGap}
The hypergraph $\HH=\HH(n,r,q|\sigma)$ with $\delta(\sigma)\geq2$ cannot have a gap above the monochromatic zone in its NMNR spectrum.
\end{proposition}

\noindent
{\bf Proof.} Let us start with $k_0$ being the largest integer for which $\HH$ has a $k_0$-NMNR colouring ($k_0$ must be less than $qn$). If $k_0=n$ or even if $k_0=n+1$ then we are done. We may therefore assume that $k_0>n+1$ and that therefore not all classes of $\HH$ are monochromatic. Choose that $k_0$-NMNR colouring of $\HH$ which has the largest number $m_0$ of monochromatic classes amongst all $k_0$-NMNR colourings.  Let $V_0$ be a non-monochromatic class of $\HH$. 

Suppose first that all colours in $V_0$ appear in some other class of $\HH$. Replace all the colours in $V_0$ with a new colour $z$. By Lemma \ref{lem:firstSwitchOfColours}, this gives an NMNR colouring of $\HH$. But this colouring uses $k_0+1$ colours. This case is therefore impossible by the maximality of $k_0$. 

Therefore suppose that $V_0$ contains just one colour $x$ which is not in some other class. Again we change all colours in $V_0$ to the colour $z$ and again we obtain an NMNR colouring, this time with $k_0$ colours but with one more monochromatic class. This is also impossible by the maximality of the number of monochromatic classes.     

Lastly, we suppose that $V_0$ contains at least two colours $x$ and $y$ which appear only in $V_0$. We replace all occurrences of these two colours in $V_0$ with the new colour $z$. By Lemma \ref{lem:secondSwitchOfColours}, we again have a NMNR colouring, but this time with $k_0-1$ colours. 

(It is important for later to observe here that the number of monochromatic classes of this $(k_0-1)$-NMNR colourings is at least $m_0$.)

We have therefore shown that $\HH$ has a $(k_0-1)$-NMNR colouring. If $k_0-1=n+1$ then we are done. 
Otherwise we again proceed as above. In general, at the $j$-th step of this procedure, we start with a $(k_0-j)$-NMNR colouring such that $k_0-j$ is still greater than $n+1$. Therefore the colouring has some non-monochromatic class $V_j$. We start with a $(k_0-j)$-NMNR colouring which has a maximum number $m_j$ of monochromatic classes. By the above observation, $m_j\geq m_{j-1}$. By re-colouring the vertices of $V_j$ as we did before, we either obtain a $(k_0-j+1)$-NMNR colouring with more than $m_j$ colourings, which is impossible since $m_j\geq m_{j-1}$, or a $(k_0-j-1)$-NMNR colouring, as required. We emphasise that the $(k_0-j-1)$-NMNR colouring we end up with has at least as many monochromatic classes as the original $(k_0-j)$-NMNR colouring, ensuring that, at the next stage, $m_{j+1}$ will be at least $m_j$.

Proceeding this way we eventually obtain that $\HH$ has a $k$-NMNR colouring with $k$ all the way from $k_0$ down to $n+1$, confirming that it has no gaps above the monochromatic zone.
\lesta

We now show, using similar techniques, that, when $\delta(\sigma)\geq2$, a $\sigma$-hypergraph cannot have gaps below the monochromatic class interval.

\begin{proposition} \label{th:noLowerGap}
Let $\HH=\HH(n,r,q|\sigma)$ be a $\sigma$-hypergraph with $\delta(\sigma)\geq2$. Then $\HH(n,r,q|\sigma)$ cannot have a gap in its NMNR spectrum below the monochromatic zone.
\end{proposition}

\noindent
{\bf Proof.} We shall proceed very much as in Proposition \ref{th:noUpperGap}. Recall that the monochromatic zone starts at $k=\lceil \frac{n}{s-1} \rceil$. Start with a $k_0$-NMNR colouring of $\HH$ with the least value of $k_0$. If $k_0= \lceil \frac{n}{s-1} \rceil -1$ then we are done. So suppose that $k_0\leq \lceil \frac{n}{s-1} \rceil -2$. Therefore the colouring has a non-monochromatic class $V_0$. Choose that $k_0$-NMNR colouring with a maximal number $m_0$ of monochromatic classes. By a suitable recolouring of the vertices of $V_0$ we either obtain a colouring with more than the maximal number of monochromatic classes or a $k_0+1$-colouring and then repeat the process. We shall describe the general case when we are at stage $j$.

In this case, we have a $(k_0+j)$-NMNR colouring where $k_0+j$ is still less than $\lceil \frac{n}{s-1} \rceil -1$. Therefore the colouring has a non-monochromatic class $V_j$. We choose a $(k_0+j)$-NMNR colouring with a maximal number $m_j$ of monochromatic classes, and we observe that $m_j>m_{j-1}$. We then have these possibilities.

If $V_j$ has at least two colours $x,y$ which do not appear in any other class of $\HH$, we colour all the vertices of $V_j$ with a new colour $z$. This gives a legitimate NMNR colouring, by Lemma \ref{lem:firstSwitchOfColours} but with strictly less colours than $k_0+j$ (say, $k_0+i$ colours, $i<j$), and one more monochromatic class (that is, $1+m_j$ monochromatic classes). This is a contradiction since $1+m_j > m_j > m_i$, for all $i<j$, therefore the new colouring has more monochromatic classes than the maximum possible for a $(k_0+i)$-NMNR colouring, which is $m_i$.

So suppose that $V_j$ has only one special colour $x$ which does not appear in any other class of $\HH$. Again we re-colour all the vertices of $V_j$ using the colour $z$, giving another $(k_0+j)$-NMNR colouring by Lemma \ref{lem:firstSwitchOfColours} but with one more monochromatic class. Again, this gives a contradiction as in the previous case.

The last remaining case is therefore when every colour in $V_j$ appears in some other class of $\HH$. We now replace all the colours in of the vertices in $V_j$ by a new colour $z$. This gives us a $(k_0+j+1)$-NMNR colouring, as required. Note that the number of monochromatic classes has also increased, insuring that, at the next step, $m_{j+1}$ will be larger than $m_j$.

Proceeding this way we finally achieve a $k$-NMNR colouring with 
$k=\lceil \frac{n}{s-1} \rceil -1$, giving us the required result. 
\lesta

From Propositions \ref{th:noUpperGap} and \ref{th:noLowerGap} we conclude the following.

\begin{theorem}
Let $\HH(n,r,q|\sigma)$ be a $\sigma$-hypergraph with $\delta(\sigma)\geq2$. Then $\HH(n,r,q|\sigma)$ cannot have a gap in its NMNR spectrum.
\end{theorem}

\section{A generalisation of NMNR colourings: colour-bounded hypergraphs} 

We have defined $\sigma$-hypergraphs and we have illustrated various cases of gaps when $\delta(\sigma)=1$. But as far as gaps go, the case for $\delta(\sigma)\geq2$ seems to be closed: there are no gaps. However, we believe that this is just the beginning of the story because  the following generalisation of NMNR colourings allows gaps to re-appear when $\delta(\sigma)\geq2$. 

In \cite{bujtas&tuz09} we find a definition which generalises NMNR colourings. A \emph{colour-bounded hypergraph} is a hypergraph together with two given integers $2\leq \alpha < \beta$ such that any feasible colouring of the hypergraph is one for which every edge has at least $\alpha$ and at most $\beta$ vertices given different colours. Since we consider only $r$-uniform hypergraphs, we also assume that $\beta\leq r$. We call such a colouring an $(\alpha,\beta)$-colouring. An $(\alpha,\beta)$-colouring using exactly $k$ colours is called a $k$-$(\alpha,\beta)$-colouring. The $(\alpha,\beta)$-spectrum is analogously defined.

Therefore a classical colouring of an $r$-uniform hypergraph is a $(2,r)$-colouring and an NMNR colouring, which is an instance of Voloshin's extension of classical hypergraph colourings, is a $(2,r-1)$-colouring. 

\subsection{The reappearance of gaps when $\delta(\sigma)\geq 2$!}

Surprisingly, just as Voloshin colourings can produce gaps when there were none in classical colourings, $(\alpha,\beta)$-colourings can produce gaps where there are none in Voloshin colourings, as the next simple result shows. 

\begin{theorem}
The hypergraph $\HH=\HH(n,4,2|\sigma=(2,2))$ with $n\geq4$ has a gap in its $(2,2)$-spectrum. 
\end{theorem}  

\noindent
{\bf Proof.} 
Colour the classes (each of size 2) of $\HH$ with the colours 1,2. Then every edge has two vertices coloured 1 and two vertices coloured 2. Therefore we have a $2$-$(2,2)$-colouring.

Next, colour the vertices in each class with the same colour, using a different colour for each class. Again, every edge contains two colours with two pairs of vertices coloured differently. This therefore gives an $n$-$(2,2)$-colouring.  

But we now show that $\HH$ does not have a $3$-$(2,2)$-colouring. So, assume that $\HH$ has been coloured with three colours. Suppose, first, that there is a non-monochromatic class. We may assume that its vertices are coloured 1 and 2, respectively. But colour 3 must also appear in some other class. But then, these two classes give us an edge with three colours, which is not allowed in a $(2,2)$-colouring. 

However, since we are using three colours and there are more than three classes, not all classes can be monochromatic. Otherwise, by the pigeonhole principle, one colour must appear in at least two  monochromatic classes, and these would give us a monochromatic edge. Therefore, since not all classes can be monochromatic, $\HH$ must have a non-monochromatic class, and we therefore obtain the same contradiction as above.
\lesta

\subsection{Further work}

Using similar agruments it can be shown that the hypergraph in the previous lemma does not have a $k$-$(2,2)$-colouring for any $3\leq k \leq n-1$. Therefore it has a gap from $k=3$ to $k=n-1$. But the point we want to make here is that while, in the world of NMNR colourings, we have a whole family of hypergraphs without gaps in their spectrum, but when considering  $(\alpha,\beta)$-colourings we see that gaps re-appear. 

Perhaps the most spactacular examples of mixed hypergraphs with broken spectra has been given by Gionfriddo \cite{gionfriddo04}. This paper presents hypergraphs which are non-colourable either for all the odd numbers or for all the even numbers between $\chi$ and $\uchi$, but colourable for all the other numbers in this range. These hypergraphs, depending on designs for their construction, are, however, not easy to produce, they are all 3-uniform, and it does not seem at all evident that it is possible to extend their construction to $r$-uniform hypergraphs with $r>3$. 

Therefore, to facilitate our study of gaps we have introduced $\sigma$-hypergraphs which are very simple to describe and which, by an appropriate choice of parameters, readily give $r$-uniform hypergraphs with $r\geq3$, with or without broken spectra. 
This led us to the discovery that $\sigma$-hypergaphs do not have gaps when $\delta(\sigma)\geq2$.  Motivated by this result and the feeling that the phenomenon of  broken spectra has not yet revealed all its secrets, we considered the $(\alpha,\beta)$-colouring of $\sigma$-graphs. With this new definition of colourings of hypergraphs it turned out that $\sigma$-hypergraphs with $\delta(\sigma)\geq2$ can have gaps. 

Our initial investigations seem to indicate that the relationship between the paramaters $\alpha$ and $\beta$ and the parameters of $\sigma$-hypergraphs can be made to work together in order to control the existence of gaps. This interplay between parameters seems to be an area of study which might give very interesting results including also a fruitful context within which one can explain why gaps in the chromatic spectra of mixed hypergraphs appear. This is work which we intend to present in a further paper.

\bibliography{colouringsAsSubmittedDM25_07_2013}

\end{document}